# Mathematik und Ethik: eine Überlegung für zukünftige Lehrpersonen

András Bátkai

*Mathematisches Modellieren wird in jedem Bereich unseres Lebens angewandt. Wird mit Hilfe von diesen Modellen über Menschenleben und Schicksale entschieden, so kommen natürlich ethische Fragen auf. Wir bieten hier eine kurzgefasste Liste von möglichen Problemen, die man auch im Schulunterricht thematisieren kann. Die angeführten Themen haben wir mit Studierenden der Mathematik in Seminaren erprobt und diskutiert.*

**Schlagwörter:** Ethik der Mathematik, Fächerübergreifender Unterricht, Studierendenprojekt

## Einleitung

Ethik wird in den österreichischen Schulen als Pflichtfach eingeführt, und zukünftige sowie praktizierende Mathematiklehrpersonen sollten auch über die ethischen Implikationen ihres Fachs nachdenken.

Wie wir das in den letzten Jahren erlebt haben, werden mathematische und statistische Argumentationen in politischen und gesellschaftlichen Diskussionen immer öfter herangezogen und teilweise auch missbraucht (z.B. Klimawandel, Corona Krise, usw.). Mathematiklehrerinnen- und Lehrer gehören zu den mathematisch gebildeten Mitgliedern unserer Gesellschaft, die eine Verantwortung tragen einerseits auf falsche Argumentationen hinzuweisen, andererseits die Wissenschaftskommunikation zu verbessern (Ernest, 2021a).

Es gab in den letzten Jahren viele Forschungsbestrebungen im Bereich Ethik der Mathematik, und bis heute sind viele spannende Fragen nicht ausreichend diskutiert worden. Unser Ziel wird es sein aufzuzeigen, welche Themen der Ethik im Mathematikunterricht angesprochen werden könnten. Dadurch wird einerseits ein Ansatz für fächerübergreifender Unterricht gegeben, andererseits wird die im Lehrplan auch geforderte gesellschaftliche Relevanz der Mathematik aufgezeigt.

Wir werden uns im Folgenden damit befassen, welche Themen man in der Lehrer*Innenbildung ansprechen sollte, damit einerseits die zukünftigen Lehrpersonen das in ihrem Unterricht einbauen können, anderseits die gesamtgesellschaftliche Verantwortung der Mathematiklehrenden bewusst wird.

Um dieses Ziel zu verwirklichen, werden wir zuerst knapp die für die Schule relevanten ethischen Themen und danach die wichtigsten Punkte der Mathematikphilosophie zusammenfassen. Danach listen wir mathematische Themen in aufsteigender philosophischer Komplexität auf, die unser Meinung nach wesentlich für die Diskussion sind. Das Thema wurde in den letzten Jahren in einem Seminar über Geschichte und Philosophie der Mathematik behandelt. Für lebhafte Diskussionen und konstruktive Anregungen bin ich meinen Studierenden der Universität Innsbruck und der PH Vorarlberg zum Dank verpflichtet. Anregungen von meinen Kolleginnen und Kollegen an der PH Vorarlberg waren auch sehr wertvoll.

## Ethik im Schulunterricht

Ethikunterricht wird in den österreichischen Schulen weiter ausgebaut, und ist somit ein wesentlicher Teil der Philosophiebildung an den Schulen. Wir fassen hier die für die Mathematik relevantesten Punkte zusammen (Pfister & Zimmermann, 2016). Ethische Fragen sind auch wichtig für die allgemeine Philosophie, weil sie metaphysische, epistemologische usw. Überlegungen gut motivieren können. Schon bei Platon fangen die Diskussionen mit Sokrates oft über ethische Fragen, wie Gerechtigkeit, an und führen zu allgemeineren Themen. Ethische Fragen haben unmittelbare gesellschaftlichen Relevanz, sowohl individuell als auch politisch. Zweck des Unterrichts ist moralische Urteile bilden zu können.

Die wichtigsten Theorien, die angesprochen werden sollten, sind die Tugendethik von Aristoteles,

die Pflichtethik von Kant, Utilitarismus, Diskursethik und Mitleidsethik. Wichtig sind die Bereichsethiken und die dazu gehörenden speziellen Fragen, wo auch die Mathematik andocken kann. Moralisch argumentieren sollte auch gelernt werden, dazu gehören logische Grundlagen wie Prämissen und Folgerungen, Syllogismus und bekannte Fehlschlüsse.

Praktische Ethik sollte problemorientiert behandelt werden. Die wichtigsten didaktischen Überlegungen sind:
- es sollte interessant und relevant sein,
- geeignetes Material für die Bearbeitung sollte zugänglich sein und
- das relevante Fachwissen sollte bei den Schülerinnen und Schülern vorhanden sein.

Unsere Auswahl von Themen versucht diesen didaktischen Überlegungen zu folgen.

### Mathematische Praktiken

Da Ethik eng damit zusammenhängt, wie Menschen handeln, sollten wir uns zuerst davon ein Bild machen, welche Arten von mathematischen Tätigkeiten es gibt. Ernest (2021a) identifiziert, in Anlehnung an die Tugendtheorie von McIntyre (Ernest, 2021b), folgende mathematischen Praktiken:
1. Rein mathematische Praxis mit Mathematik explizit als Ziel. Hierhin gehört theoretisch-mathematische Forschung.
2. Angewandte mathematische Praxis mit Mathematik explizit als Ziel. Beispiele wären Forschung in angewandter Mathematik.
3. Angewandte und theoretische Praxis mit Mathematik implizit. Mathematisch basierte theoretische Forschung in anderen Wissenschaften wäre so ein Beispiel, wie theoretische Physik.
4. Rein praktische Aktivität mit Mathematik explizit als Ziel. Unterricht und die Entwicklung von Unterrichtsmaterialien gehören beispielsweise hierhin. Aber auch Mathematiklernen hat Mathematik explizit als Ziel.
5. Reine und angewandte theoretische Arbeit, Mathematik implizit. Die Entwicklung von Apps und Computeranwendungen, die in unterschiedlichsten Lebensbereichen eingesetzt werden, wäre ein Beispiel für diesen Punkt.
6. Angewandte und praktische Tätigkeiten, Mathematik explizit. Wirtschaftsmathematik, Finanzstatistik, usw. sind Beispiele hierfür. Die Tätigkeit konzertiert sich auf die Entwicklung neuer Anwendungen.
7. Reine und angewandte praktische Arbeit, Mathematik implizit. Dazu gehören Alltagsmathematik wie Steuern, Gehaltsberechnung, usw. Allgemeiner, Datenverarbeitung und Verwendung in unserem täglichen Leben könnte man auch hier platzieren.

Die einzelnen Punkte sind hier diskutierbar, wie zum Beispiel ob Punkte 1 und 2 nicht zusammengehören, aber die Liste gibt eine ziemlich gute Zusammenfassung darüber wie und wo heute Mathematik betrieben wird. Wir können uns fragen, welche dieser Punkte in unserem Unterricht thematisiert werden, wenn wir Mathematik als lebendige Wissenschaft darstellen, wie der Lehrplan das vorschreibt. Im Folgenden werden wir uns auf die ethischen Aspekte der mathematischen Modellierung konzentrieren, und hier die für die praktische Ethik relevantesten Punkte, also 5 bis 7, detaillierter thematisieren. Wir werden dabei stets den schulischen Kontext im Auge behalten.

### Themen

Im Folgenden werden wir, mit vielen Beispielen, die im Unterricht behandelt werden können, solche mathematische Themenfelder betrachten, die ethisch relevante Fragen aufwerfen. Wir machen hier keine Taxonomie, einige Beispiele könnten durchaus in mehreren Kategorien erscheinen, und die einzelnen mathematische Themen sind auch nicht voneinander getrennt. Die Problemfelder sind nach didaktischen Überlegungen in aufsteigender philosophischer Komplexität präsentiert.

#### Rechenfehler
Die offensichtlichste Situation ist, wenn wir Beispiele dafür anschauen, was für Folgerungen einfache Rechenfehler haben können. Es ist wichtig, uns bewusst zu machen, dass bei den Anwendungen der Mathematik Rechenfehler zu verheerenden

Situationen führen können. Obwohl Fehler zur Schulalltag dazugehören und einen wichtigen didaktischen Beitrag leisten, beim entdeckenden Lernen sogar eine zentrale Funktion haben, sollten wir diesen Aspekt im Unterricht nicht vergessen.

- Der Mars climate orbiter (Grossman, 2010) ist im Jahr 1999 kaputtgegangen, weil bei der Vorbereitung ein Zulieferer die Maße nicht in SI-Einheiten gerechnet hat, wodurch die Teile nicht kompatibel waren.
- Vasa (Hocker, 2011) war ein Kriegsschiff der königlichen Marine, die bei seiner ersten Ausfahrt im Jahr 1628 Stabilitätsprobleme hatte und untergegangen ist. Bei der Untersuchung wurde festgestellt, dass mehrere Handwerker in unterschiedlichen Maßeinheiten gerechnet (Schwedisches Fuß und Amsterdamer Fuß) haben und so Teile des Schiffes nicht wirklich zusammengepasst haben. Bei dem Untergang sind um die 50 Menschen getötet worden.
- Die Hochrheinbrücke (Der Spiegel, 2004) wurde von zwei Seiten aus gebaut um als Verbindung zwischen Deutschland und der Schweiz zu dienen. Durch einen Vorzeichenfehler hatten aber die Teile 56 cm Höhenunterschied und die Brücke musste zurückgebaut werden. Dabei sind Kosten in Millionenhöhe entstanden.

### Plausibilität

Viele Entscheidungen des täglichen Lebens hängen davon ab, wie plausibel bestimmte mögliche zukünftige oder vergangene Ereignisse uns erscheinen. Die Schuldfrage vor Gericht wird meistens danach entschieden, wie glaubhaft es dargestellt wird, dass der Angeklagte die Tat wirklich begangen hat (Buchanan, 2007). Da es hier um Gerechtigkeit und menschliche Schicksale geht, ist die Beziehung zur Ethik klar gegeben. Um mehr Objektivität zu schaffen, ist die Versuchung hier groß Plausibilität mit (mathematischen) Wahrscheinlichkeit zu vertauschen, zu modellieren. Also ein Ereignis wird als glaubhaft angesehen, wenn derer mathematische Wahrscheinlichkeit groß ist. Die Alltagsprache suggeriert diesen Schritt auch: wir benutzen das Wort „wahrscheinlich" als synonym für „plausibel". Pólya (1975) hat in seinem berühmten didaktischen Werk auch versucht, Plausibilität mit Wahrscheinlichkeit zu erklären, und der Neymar-Pearson Test in der Statistik basiert auch darauf, dass man Plausibilität durch Wahrscheinlichkeit ersetzt.

Arbeitet man mit Wahrscheinlichkeiten, kann es vorkommen, dass es zu fehlerhaften Argumentationen kommt. In den folgenden Beispielen wird die fehlerhafte Verwendung von bedingten Wahrscheinlichkeiten thematisiert, und die meisten demonstrieren das, was als „Trugschluss des Anklägers" oder „Konfusion des Staatsanwalts" (prosecutor's fallacy (Krämer, 2015) in die Literatur eingegangen ist. Siehe auch (Schneps & Colmez, 2013).

- Sally Clark (Byard, 2004) wurde 1999 angeklagt ihre zwei neugeborenen Kinder umgebracht zu haben, wo die ärztliche Diagnose zuerst plötzliche Kindestod lautete. Der Hauptzeuge der Anklage war ein Arzt, dessen Argument wie folgt lautete: Die Wahrscheinlichkeit für ein plötzliches Kindestod ist 1:8500, als für zwei Kinder ist sie das Quadrat. Also ist es sehr unwahrscheinlich, dass beide Kinder an plötzlichen Kindestod starben, so muss Frau Clark ihre Kinder ermordet haben. Sally Clark wurde verurteilt und es hat drei Jahre gebraucht ihre Unschuld in Revisionsverfahren zu beweisen. Sie hat sich leider nie wieder erholt und ist später an Alkoholvergiftung verstorben. Ihr Fall ist ein besonders schockierender, der zeigt, wie viel Schaden man mit falschen mathematischen Argumenten und fehlerhaften Wahrscheinlichkeitsberechnungen anrichten kann.
- Lucia de Berk (Buchanan, 2007) hat als Krankenschwester gearbeitet und wurde verdächtigt, einige ihre Patienten getötet zu haben. Erstaunlich am ganzen Verfahren war, dass nie auch nur ein konkreter Beweis für ein Verschulden Lucia de Berks aufgeführt wurde. Unter anderem wurden dafür Wahrscheinlichkeitsrechnungen angestellt, die nachweisen sollten, dass eine zufällige Anwesenheit der Pflegerin bei den verdächtigen Todesfällen eigentlich unmöglich sei. Die Argumente waren statistischer Natur, und die Statistischen Fehler, die man im Prozess begangen hat, wurden in der Arbeit (Meester et al., 2007) sehr gut, auch teilweise für die Schule benutzbar, aufgearbeitet. Eines von den mehreren falschen Argumenten war, dass es sehr unwahrscheinlich war, dass sie bei all den Todesfällen nur zufällig anwesend gewesen war, also muss sie die Tode verursacht haben. Erst nach langjährigen Gerichtsverfahren wurde sie freigesprochen.
- O.J. Simpson wurde freigesprochen (Krämer, 2015), seine Frau ermordet zu haben, obwohl die Öffentlichkeit von seiner Schuld überzeugt war. Es war nicht wesentlich im Prozess, aber seine Verteidigung hat unter Anderem falsche Argumente gebracht und hat versucht, zu argumentieren, dass nur weil Simpson seine Frau geschlagen hat, er sie nicht ermordet haben muss, denn viele Männer schlagen ihre Frauen aber nur wenige Töten sie. Hätte man aber die Frage gestellt, bei wie vielen getöteten Frauen, die von ihrem Mann geschlagen worden sind, auch von ihm getötet sind, wäre ein ganz anderes Bild entstanden.

### Falsche Modelle oder Fehlinterpretation der Modelle

Dieses Thema ist meistens nicht so gut in der Schule zu behandeln, weil die Modelle meist kompliziert und mathematisch anspruchsvoll sind. Es ist trotzdem wichtig, wenn wir uns bewusst machen, dass solche Phänomene in den öffentlichen Diskussionen häufiger vorkommen.
- Die Finanzkrise im Jahr 2009 hat die Welt erschüttert, und viele machten die (falsch angewandten) mathematischen Modelle verantwortlich. Die Literatur, in der die ethischen Implikationen aufgearbeitet werden, ist umfangreich, siehe z.B. Ippoliti (2021) oder Johnson, (2012).
- Manipulation mit Hilfe von Statistiken ist ein wichtiges Thema, siehe (Krämer, 2015), der uns bei politischen Diskussionen regelmäßig begleitet. Man denke nur an die verschiedenen Interpretationen von unterschiedlichen Kennzahlen wehrend der Corona-Pandemie.
- Die Corona-Pandemie bietet auch zahlreiche Beispiele, wo mathematische Vorhersagen fehlinterpretiert worden sind, denkt man nur an die Inzidenz, die Sicherheit der Tests oder der Wirkungsgrad der Impfungen. Dass die mathematischen Modelle nicht gut funktioniert haben (Ioannidis et al., 2022), weil wir zu wenig über das Virus gewusst haben, hat auch zur Unsicherheit beigetragen.

### Textaufgaben
- Solange Schulaufgaben „rein" mathematisch sind, kann man davon ausgehen, dass keine großen ethischen Diskussionen durch die Aufgabentexte angestoßen werden. Im Unterricht werden aber auch viele Aufgaben mit Alltagsrelevanz behandelt, und solche Textaufgaben sind ethisch nie neutral, siehe (Shulman, 2002). Oft vermitteln sie ein Wertesystem und suggerieren Haltungen. Als Illustration wird hier die Extremsituation des Nationalsozialismus behandelt, wo man in Schulbüchern reichlich Beispiele (Heske, 2021) dafür findet, in denen berechnet werden soll,
- Was der Anteil der Juden in unterschiedlichen Bevölkerungsschichten ist.
- Wie viele Familien könnten Häuser kaufen aus dem Geld, das der Staat für Geisteskranke ausgibt.
- Wie viele Leute in einem Luftschutzkeller, wie viel Bomben in einem Flugzeug passen.
- Welche Prozente bei der Saarabstimmung für Deutschland gestimmt haben.

Die Lehrpersonen sollten sich aber bei jeder gestellten Aufgabe Fragen, welche Werte vermittelt werden und ob eine ethische Diskussion der Aufgabe auch angebracht ist. Kontexte wie Klimawandel, Umwelt, Gesundheit und Armut laden sogar dazu ein.

### Fragwürdige Anwendungen
Mathematik wird seit seinen Anfängen für militärische Zwecke genutzt und durch militärische Anwendungen motiviert. Man findet schon in altbabylonischen Keilschrifttexten Aufgaben, wo man das Volumen von Belagerungsdämmen berechnen soll, und in jede Epoche der Geschichte der Mathematik ist die enge Beziehung zu militärischen Anwendungen mehr als nur sichtbar, man denke nur an Archimedes. Das Buch (Booß-Bavnbek & Hoyrup, 2003) gibt einen guten Überblick, wie man überall Mathematik im Krieg verwendet hat.

### Automatisierte Entscheidungen
Soll man zwischen verschiedenen Optionen wählen, wünschen sich viele Algorithmen, die diese Wahl automatisch treffen. Man wünscht sich dadurch mehr Objektivität und Gerechtigkeit. Eines der Probleme mit automatisierten Entscheidungen ist, dass entweder die Algorithmen nicht öffentlich sind (Betriebsgeheimnis), und somit die Transparenz der Entscheidungen nicht gewährleistet ist, oder, wenn sie öffentlich gemacht werden, entsteht oft eine Rückkopplungsschleife und somit eine Verfälschung des Systems. Es wird auch oft beobachtet, dass automatisierte Entscheidungen Ungleichheiten sowie Stereotypien verstärken, da sie durch Datenverarbeitung nicht Kausalität, sondern Korrelationen feststellen.

- Algorithmen kommen immer öfter zum Einsatz, wenn bei in-Vitro Fertilisation Embryonen ausgewählt werden (Afnan et al., 2021). Das Thema ist mathematisch schwierig und ethisch fragwürdig. Meistens werden Wahrscheinlichkeiten berechnet und Embryonen mit den besten Überlebenschancen ausgewählt. Dieser utilitaristische Ansatz wirft mehr Fragen als Antworten auf.
- Um Verbrechen vorzubeugen, werden immer öfter auch Algorithmen herangezogen (Bringsjord et al., 2021). Hier werden meistens durch Big-Data Analysis auch Wahrscheinlichkeiten berechnet. Und wenn zum Beispiel in einer Gegend bis jetzt Einbrüche öfter von Mitgliedern einer ethnischen Gruppe begangen worden sind, dann wird schon die Zugehörigkeit zu dieser Gruppe als Risikofaktor bewertet und so die Vorurteile bestärkt.
- Autonomes Fahren hat im Bereich der Technikethik seinen besonderen Stand und auch reichlich Literatur. Wie Autos in Gefahrensituationen reagieren sollen, ist noch weit weg davon ausdiskutiert worden zu sein, siehe (Funk, 2022).
- Ein für den Unterrichtsbereich heikles Thema sind standardisierte Tests. Die bieten zwar Vergleichbarkeit und Objektivität, haben aber ein Rückkopplungseffekt auf den Unterricht, was nicht immer willkommen ist. Es gibt viele Berichte zum Beispiel dazu, dass Themen, die zwar im Lehrplan stehen, aber in der Zentralmatura wahrscheinlich nicht geprüft werden, kaum im Unterricht thematisiert werden. Das Spannungsfeld zwischen Standardisierung und Individualisierung sollte auch nicht außer Acht gelassen werden (Brügelmann, 2019).
- Suchmaschinenalgorithmen werden gerne als schöne mathematische Beispiele, zum Beispiel bei linearer Algebra angesehen, werden aber auch heftig kritisiert wegen deren Beitrag zur Verbreitung von Verschwörungstheorien (Ghezzi et al., 2020). Es werden auch Stereotypien, ob bestimmte Berufe lieber von Männern oder von Frauen ausgeübt werden, bestärkt (Kay et al., 2015), und generell, die Manipulierbarkeit der Suchergebnisse (Epstein & Robertson, 2015) ist ein großes Problem.
- Allgemeiner, Ranking hat seinen eigenen Bereich reichlich belegt mit Literatur, wenn es um die Diskussion der ethischen Auswirkungen geht. Universitätsrankings, Arbeitgeberrankings usw. haben ihre eigene Problemfelder (O'Neil, 2017).

**Modellierung**

Die Mathematisierung und mathematische Modellierung von Phänomenen um uns herum hat eine Geschichte fast gleich so alt wie die menschliche Kultur. Die philosophische Grundposition, die dahintersteckt, nämlich, dass unsere Welt (approximative) erkennbar ist (Nickel, 2007), möchte ich hier nicht diskutieren, aber der Modellierungsprozess wirft ethische Fragen auf, und die sollten im Unterricht auch thematisiert werden. Falsche Vorstellungen über mathematische Modelle tragen sicherlich zur Wissenschaftskepsis bei, auch wenn das Problem viel komplexer ist (Kollosche, 2021).

Den Modellierungsprozess möchte ich durch die Ausführungenl von Andreas Vohns (2010) thematisieren. Das ist ein mathematisch einfaches Beispiel, wo die Rolle der Mathematik sich auf die Verwendung von Zahlen beschränkt. Aber genau aus diesem Grund kann man die einzelnen Modellierungsschritte besser erkennen.

- Man möchte in einem Wohngebiet vorsichtiges und verantwortungsvolles Fahren modellieren, mit dem Ziel, die Autofahrer zu solchem Verhalten zu motivieren und verantwortungslose Autofahrer zu bestrafen.
- In unserem Kulturkreis ist es üblich, das Vorherige durch Geschwindigkeit zu modellieren. Wenn man gar nicht fährt, ist es die absolute Fahrsicherheit, und wenn man rast, dann ist man verantwortungslos.

- Jetzt kommt die Frage nach der Geschwindigkeitsbegrenzung. Wenn man die Autos durchfahren lässt, dann ist 0 keine Option, weil man weiterkommen möchte. Also die Begrenzung wird ein Kompromiss zwischen Sicherheit und Vorwärtskommen. Es gibt zwar viele Argumente für bestimmte Begrenzungen, aber keine von Natur aus zwingende. Man kann 20, 30 oder 40 oder was anderes argumentieren.
- Setzt man eine Zahl, zum Beispiel 30, fest, wird diese für alle zwingend. Fährt jemand schneller, der kann bestraft werden und die Strafe muss man auch dann bezahlen, wenn man eine andere Grenze besser argumentieren kann.

Also geschieht Modellierung in der Dialektik zwischen Willkür und Zwang. Mathematische Modellierung wird eingesetzt, weil die Ergebnisse objektiv nachvollziehbar sind, und weil man ein hohes Maß an Disziplinierung dadurch erreichen kann.

- Bei der Corona-Pandemie hat man auch versucht, die Auslastung der Intensivstationen zu modellieren. Dabei hat man (basierend auf Erfahrungen mit der Spanischen Grippe) zuerst mit der Inzidenz modelliert, siehe Cao & Liu (2022). So hat man versucht zwischen dem sicheren Schutz, d.h. vollständiger Lockdown und völliger Freiheit eine Linie zu ziehen, die sowohl ein gewisses Maß an Sicherheit als auch funktionierenden Wirtschaft bietet. Da dieser Konsens in unterschiedlichen Ländern zu unterschiedlichen Maßnahmen geführt hat, konnte die Politik den zwingenden Teil der Modellierung nicht immer gut kommunizieren.
- Menschenwürde soll in den meisten west-europäischen Ländern geachtet und geschützt werden. Die meisten in unserem Kulturkreis sind auch davon überzeugt, dass der Staat dabei eine Fürsorgefunktion hat und die bedürftigen unterstützen soll. Bei der Modellierung kommt man schnell auf die Idee, menschenwürdiges Leben durch Einkommen zu charakterisieren und die staatliche Hilfe wird als Sozialhilfe realisiert, siehe Vohns (2010). Dabei ist die Spanne zwischen keine Sozialhilfe und ein Grundeinkommen für alle groß. Abhängig von dem politischen Diskurs kann man ein Warenkorbmodell oder ein Statistikmodell wählen, um in einem Diskurs menschenwürdiges Existenzminimum zu definieren.
- Oft ist es so, dass Entscheidungsträger eindimensionale Schwellenwerte für bestimmte Kennzahlen brauchen, damit Entscheidungen leicht begründet und nicht mehr diskutiert werden können. Dabei sollte man die willkürlichen Entscheidungen im Modellierungsprozess nicht aus den Augen verlieren und Raum für individuelle Ausnahmeentscheidungen haben.

### Fazit

Kant wird oft zitiert: „Ich behaupte aber, dass in jeder besonderen Naturlehre nur so viel eigentliche Wissenschaft angetroffen werden könne, als darin Mathematik anzutreffen ist" (Kant, 1786). Es scheint so, als ob dieser Satz sein Eigenleben begonnen hätte und wir erleben eine noch nie dagewesene Mathematisierung der Welt (Nickel, 2007; 2022), die schon lange die ursprünglichen Intentionen von Kant hinter sich hat. Jeder Lebensbereich von uns wird mit Zahlen charakterisiert in der Hoffnung, Objektivität und Gerechtigkeit zu schaffen. Wie wir aber an den Beispielen gesehen haben, erreicht man oft das Gegenteil: Unterschiede werden verstärkt und die Ungerechtigkeit vergrößert. Aus diesem Grund ist es vital, dass zukünftige Mathematiklehrpersonen sich mit den ethischen Fragen mathematischer Modellierung beschäftigen. Bei unserem Seminar haben die Studierenden sich in diese und viele andere Themen, wie Garrymandering, Triage, uvm. vertieft und bei Workshops kritisch diskutiert. Sie haben sich überlegt, wo im Unterricht die Themen reinpassen könnten und haben auch Unterrichtssequenzen erstellt.

Andere Beziehungen zwischen Mathematik und Ethik, nämlich Professionsethik (Plagiat, Priorität, Minoritäten, uvm.), oder die Fragestellung, was eigentlich gute (und schlechte) Mathematik ist, konnten wir hier nicht behandeln. Die Letztere, die sehr eng mit der Frage, was Mathematik ist, zusammenhängt, spielt für den Schulunterricht eine wesentliche Rolle. Sie entscheidet darüber, was wir unterrichten und wie wir unterrichten (Ernest, 2018; Maaß & Götz, 2022). Das sind aber Themen für eine andere Veröffentlichung.